\documentclass[11pt, reqno]{amsart}

\usepackage{amsmath,amssymb,mathtools}
\usepackage{bm}
\usepackage[colorlinks=true,linkcolor=blue,citecolor=blue]{hyperref}
\makeatletter
\let\mark@size\scriptsize
\makeatother
\usepackage{enumitem}
\usepackage{xcolor}

\usepackage[margin=1.15in]{geometry}

\theoremstyle{plain}
\numberwithin{equation}{section}

\newtheorem{theorem}{Theorem}[section]
\newtheorem{lemma}{Lemma}[section]
\newtheorem{corollary}{Corollary}[section]
\newtheorem{proposition}{Proposition}[section]
\newtheorem{remark}{Remark}[section]

\newtheorem*{notation}{Notation}

\newcommand{\T}{\mathbb{T}}
\newcommand{\R}{\mathbb{R}}
\newcommand{\N}{\mathbb{N}}

\newcommand{\pa}{\partial}
\newcommand{\grad}{\nabla}

\newcommand{\Dsym}{D}
\newcommand{\nrm}[2]{\|#1\|_{#2}}

\newcommand{\abs}[1]{|#1|}

\newcommand{\intt}{\int_{\T^3}}
\newcommand{\Lp}[1]{L^{#1}}

\newcommand{\calC}{\mathcal C}

\newcommand{\calE}{\mathcal{E}}
\newcommand{\calD}{\mathcal{D}}

\newcommand{\e}{\varepsilon}

\newcommand{\Gv}{G[\abs{\Dsym v}^2]}

\newcommand{\lt}{\left}
\newcommand{\rt}{\right}
\newcommand{\bq}{\begin{equation}}
\newcommand{\eq}{\end{equation}}
\newcommand{\sfI}{\mathsf{I}}

\begin{document}

\title[pressureless Euler/non-Newtonian on $\T^3$]{Global-in-time dynamics of the 
pressureless Euler/incompressible non-Newtonian equations}

\author[Jung]{Jinwook Jung}
\address[Jinwook Jung]{\newline Department of Mathematics and  Research Institute for Natural Sciences \newline
Hanyang University, 222 Wangsimni-ro, Seongdong-gu, Seoul 04763, Republic of Korea}
\email{jinwookjung@hanyang.ac.kr}

\author[Kim]{Subin Kim}
\address[Subin Kim]{\newline Department of Mathematics \newline
Hanyang University, 222 Wangsimni-ro, Seongdong-gu, Seoul 04763, Republic of Korea}
\email{ksb6575@hanyang.ac.kr}

\allowdisplaybreaks

\date{\today}


\begin{abstract}
This paper addresses a two-phase fluid system on the periodic domain $\T^3$, consisting of the pressureless Euler equations and the incompressible non-Newtonian equations coupled through a drag force. We prove the global-in-time existence of strong solutions to this system for a certain class of stress tensors in the non-Newtonian fluid. Our proof combines a priori large-time behavior estimates for classical solutions with a bootstrapping argument on the $L^2$-Sobolev norms, which effectively handles the difficulties caused by the lack of dissipation in the density of the Euler equation.
\end{abstract}

\maketitle

\tableofcontents

\section{Introduction}
\setcounter{equation}{0}

In this paper, we consider the global-in-time existence and uniqueness of strong solutions to a coupled hydrodynamic system on the periodic domain $\T^3=\R^3/\mathbb{Z}^3$. The system of our concern comprises the pressureless Euler equations coupled with incompressible, non-Newtonian fluid equations via a drag force. More specifically, let $\rho= \rho(t,x)$ and $u=u(t,x)$ be the density and velocity for the pressureless flow at $(t,x) \in \R_+\times \T^3$, respectively, and $v=v(t,x)$ be the velocity for the incompressible flow. Then the triplet $(\rho, u, v)$ solves the following system:

\bq\label{main_sys}
\begin{aligned}
&\pa_t \rho +\nabla_x \cdot (\rho u)=0, \quad (t,x)\in \R_+ \times \T^3,\\
&\pa_t (\rho u)+\nabla_x \cdot (\rho u\otimes u)=\rho(v-u),\\
&\pa_t v +(v\cdot\nabla_x)v+\nabla_x p - \nabla_x\cdot \!\bigl(\Gv\,\Dsym v\bigr)=-\rho(v-u),\\
&\nabla\cdot v=0,
\end{aligned}
\eq
subject to initial data:
\bq\label{init_data}
(\rho(0,x), u(0,x), v(0,x)) =: (\rho_0(x), u_0(x), v_0(x)), \quad x \in \T^3.
\eq
Here, $p=p(t,x)$ and $G = G[|Dv|^2]$ denote the pressure and the viscous part of the stress tensor of the incompressible fluid, respectively, and $Dv$ denotes the symmetric part of the velocity gradient, i.e. $Dv = (\nabla v + \nabla v ^T)/2$. For the viscous part $G= G(s) : \R_+ \to \R_+$ of the stress tensor, we assume that $G$ is a smooth function and impose the following structural conditions: for any $s \in [0,\infty)$,

\begin{align}
\tag{H1}\label{H1} & G[s]\ge m_0>0 \\
\tag{H2}\label{H2} & G[s]+2G'[s]\,s\ge m_0>0 \\
\tag{H3}\label{H3} & \bigl|G^{(k)}[s]\,s^{\alpha}\bigr|\le C_k\,\bigl|G^{(k-1)}[s]\bigr|,
\quad \alpha\in\{0,1\},\ k\in\N 
\end{align}
where $G^{(k)}[\cdot]$ is the $k$-th order derivative of $G$ and $m_0$ and $C_k$ are positive constants.
A typical example for $G$ satisfying \eqref{H1}-\eqref{H3} is the following power-law type:

$$G[s]=(m_0^{2/(q-2)}+s)^{(q-2)/2} \quad \mbox{with } \ 2<q<\infty. $$

System \eqref{main_sys} describes the dynamics of particles dispersed through a (non-Newtonian) fluid and is closely related to kinetic-fluid coupled systems which have drawn substantial interest, motivated by modelling sedimentation phenomena, spray dynamics, combustion, etc. \cite{BWC00, R81, W58}. In our case, one can find a connection between system \eqref{main_sys} and the Vlasov/non-Newtonian coupled system: 

\bq\label{eq:kf_sys}
\begin{aligned}
&\pa_t f + \xi \cdot \nabla_x f +\nabla_\xi\cdot ((v-\xi) f) = 0,\quad (x,\xi)\in \T^3 \times \R^3, \quad t>0,\\
&\pa_t v + (v\cdot \nabla_x) v + \nabla_x p - \nabla_x\cdot \!\bigl(\Gv\,\Dsym v\bigr)=-\int_{\R^3}(v - \xi)f \,d\xi,\\
&\nabla\cdot v = 0,
\end{aligned}
\eq
where $f = f(t,x,\xi)$ denotes the probability density of dispersed particles in phase space $(x,\xi)\in \T^3 \times \R^3$ at time $t>0$ and $v = v(t,x)$ denotes the velocity of the incompressible flow. Once we write the zeroth- and first-order velocity moments as
\[
\rho_f(t,x) :=\int_{\R^3} f(t,x,\xi)\,d\xi, \quad (\rho_f u_f )(t,x) := \int_{\R^3} \xi f(t,x,\xi)\,d\xi,
\]
we deduce from \eqref{eq:kf_sys} that $\rho_f$ and $\rho_f u_f$ satisfy, at the formal level,
\[
\begin{aligned}
&\pa_t \rho_f + \nabla_x \cdot (\rho_f u_f) =0,\\
&\pa_t (\rho_f u_f) + \nabla_x \cdot (\rho_f u_f \otimes u_f) + \nabla_x \cdot \lt(\int_{\R^3} (\xi - u_f)\otimes (\xi-u_f) f \,d\xi \rt) = - \rho_f (u_f - v),\\
&\pa_t v + (v\cdot \nabla_x) v + \nabla_x p - \nabla_x\cdot \!\bigl(\Gv\,\Dsym v\bigr)= -\rho_f (v - u_f),\\
&\nabla\cdot v = 0,
\end{aligned}
\]
By assuming the monokinetic ansatz 
\[
f(t,x,\xi) = \rho_f(t,x)\,dx\otimes \delta_{\xi = u_f(t,x)}(d\xi),
\]
the equations for $\rho_f$ and $\rho_f u_f$ formally reduce to the pressureless Euler equations. A rigorous derivation of \eqref{main_sys}, however, remains largely open due to the lack of sufficient control on $\rho_f$ even in the case of a linear viscous stress tensor $G \equiv m_0 >0$. On the other hand, several works have addressed hydrodynamic limits from kinetic-fluid models to two-phase fluid systems; see \cite{CCK16,CJ21-1,MV08} and the references therein.

Independently of such rigorous derivations, fluid-fluid coupled systems have been extensively studied. 
In these models, the dispersed particles are typically described by compressible Euler equations, while the surrounding fluid is governed by incompressible or compressible Navier-Stokes equations. In \cite{Choi15, Choi16, HTWZ24}, compressible Euler equations coupled to incompressible/compressible Navier-Stokes equations were considered and global existence of strong solutions was established, together with the large-time behavior of velocity fields. In this case, the pressure term in the Euler equation provides proper control on the Euler density from coercivity-type estimates. Combined with $L^2$-type energy estimates, this mechanism plays a crucial role in the construction of the global-in-time solutions.

In the absence of pressure, however, no analogous dissipation estimates for the Euler density are available. To overcome this difficulty, {\it a priori} large-time behavior estimates have been employed for the global existence of pressureless Euler equations coupled to compressible Navier-Stokes equations on the periodic domain $\T^3$ \cite{CK16} and to incompressible Navier-Stokes equations on the whole space $\R^3$ \cite{CJ21-2}, respectively.  The whole-space result was subsequently improved and extended to the case in which the incompressible Navier-Stokes equations are replaced by compressible ones; see \cite{CJK24, Dan26, GWZ24, HTZ24, Lpre, LTZ25, ZCLZ23}.

Motivated by these developments, we study a coupled fluid-fluid system in which the dispersed particles are described by the pressureless Euler equations and the surrounding fluid by the incompressible non-Newtonian fluid equations. This setting combines the difficulties associated with pressureless flow with those arising from a nonlinear viscous stress tensor. 

There has been some progress in the study of kinetic equations coupled with non-Newtonian fluid equations \cite{CJW27, HKKP18, MPP18, Wpre}. Moreover, the mathematical theory of non-Newtonian fluid equations themselves has a long and substantial history; see, for example, \cite{BW15,DR05,DRW10,KKK19,L69,MNRR96}. Nevertheless, to the best of our knowledge, the pressureless Euler equations coupled to an incompressible non-Newtonian fluid have not previously been studied. In view of the extensive literature on both coupled fluid-fluid systems and non-Newtonian fluid equations, it is quite natural to ask whether the global well-posedness framework developed for Euler/Navier-Stokes systems can be extended to the Euler/non-Newtonian setting. The purpose of this work is to give a positive answer to this question by establishing global-in-time existence and uniqueness of strong solutions to \eqref{main_sys} under suitable assumptions on the initial data.

Before we state our main theorem, we present the averaged quantities:

\bq\label{aver_quant}
\rho_c (t) := \intt \rho\,dx, \quad u_c(t) := \intt u\,dx,  \quad v_c(t) := \intt v\,dx.
\eq

Then the main theorem is stated below:

\begin{theorem}\label{main_thm}
Suppose that the initial data $(\rho_0, u_0, v_0)$ satisfy
\begin{enumerate}
\item $\rho_0(x) >0$ for every $x\in \T^3$ and $\intt\rho_0\,dx=1$, and
\item $(\rho_0, u_0, v_0) \in H^2(\T^3)\times H^3(\T^3)\times H^3(\T^3)$ with $\nabla\cdot v_0=0$.
\end{enumerate}
If
\[
\|\rho_0 -\rho_c(0)\|_{H^2} + \|u_0 - u_c(0)\|_{H^3} + \|v_0 - v_c(0)\|_{H^3} + |u_c(0)-v_c(0)| < \e_0
\]
for $\e_0>0$ sufficiently small, then system \eqref{main_sys}-\eqref{init_data} admits a unique global-in-time strong solution $(\rho, u, v) \in \calC([0,\infty); H^2(\T^3)) \times \calC([0,\infty); H^3(\T^3)) \times \calC([0,\infty); H^3(\T^3))$ satisfying $\rho(t,x) >0$ for every $(t,x) \in \R_+\times\T^3$ and
\[
\sup_{t \ge 0}\lt( \|\rho(t, \cdot) -\rho_c(t)\|_{H^2} + \|u(t,\cdot) - u_c(t)\|_{H^3} + \|v(t,\cdot) - v_c(t)\|_{H^3}\rt)<\infty.
\]
Moreover, there exist constants $C, c>0$ independent of $t$ such that
\[
\|u(t,\cdot) - u_c(t)\|_{H^3} + \|v(t,\cdot) - v_c(t)\|_{H^3} \le Ce^{-ct}.
\]

\end{theorem}

\begin{notation}
For a multi-index $\alpha$ with $|\alpha|=l$ we write $\abs{\grad^l v}^2=\sum_{\abs\alpha=l}\abs{\pa^\alpha v}^2$ and $\nrm{\grad^l v}{\Lp2}^2=\intt\abs{\grad^l v}^2\,dx$; the object $\pa^\alpha\Dsym v$ with $\abs\alpha=l$ is of order $l+1$ in $v$, so that $\nrm{\grad^l \Dsym v}{\Lp2}$ controls the $(l+1)$-st order of $v$. For $3\times 3$ matrices $A = (a_{ij})$ and $B=(b_{ij})$, we write $A:B := \sum_{i,j=1}^3 a_{ij} b_{ij}$.  For simplicity, we omit $x$-dependence of differential operators, i.e. $\nabla f := \nabla_x f$, $\Delta f := \Delta_x f$ and $\nabla\cdot F := \nabla_x \cdot F$. The symbol $C$ denotes a generic positive constant, independent of $t$, which may change from line to line.
\end{notation}

\section{Preliminaries}
\setcounter{equation}{0}

In this section, we present several technical lemmas used throughout this paper. We also provide local-in-time existence of strong solutions to \eqref{main_sys}-\eqref{init_data}.

We first recall Korn's inequality and Moser-type estimates:

\begin{lemma}
Let $w$ be a divergence-free, mean-zero smooth vector field on $\T^3$ and $f, g \in H^k(\T^3) \cap L^\infty(\T^3)$ for some $k\in \N$.
\begin{enumerate}
\item
\emph{(Korn's inequality)}
There exists a positive constant $c_K$ independent of $w$ such that

\[
\|\nabla w\|_{L^2} \le c_K \|Dw\|_{L^2},
\]
 
\item
There exists a positive constant $C>0$ independent of $f$ and $g$ such that
\[
\|\nabla^k (fg) \|_{L^2} \le C\lt( \|\nabla^k f\|_{L^2}\|g\|_{L^\infty} + \|f\|_{L^\infty}\|\nabla^k g\|_{L^2}\rt).
\]
Moreover, if $\nabla f \in L^\infty(\T^3)$, we have
\[
\|\nabla^k (fg) - f\nabla^k g\|_{L^2}\le C\lt( \|\nabla^k f\|_{L^2}\|g\|_{L^\infty} + \|\nabla f\|_{L^\infty}\|\nabla^{k-1}g\|_{L^2}\rt).
\]
 
\end{enumerate}
\end{lemma}
Next, we recall the higher-order derivative estimates and coercivity estimates for the viscous stress tensor whose proofs can be found in \cite{KKK19}.

\begin{lemma}\cite[Lemma 2.1]{KKK19}\label{lem:kkk1}
For a multi-index $\alpha$ with $\abs\alpha=\ell$,
\begin{equation}
\pa^\alpha\!\bigl(\Gv\,\Dsym v\bigr)
=\Gv\,\pa^\alpha\Dsym v
+2G'[\abs{\Dsym v}^2]\,(\Dsym v:\pa^\alpha\Dsym v)\,\Dsym v
+E_\alpha,
\end{equation}
where $E_\alpha$ collects all the lower-order terms. For $\ell=2$ and $\ell=3$, $E_\alpha$ also obeys the following pointwise bounds:
\bq
\begin{aligned}
|E_\alpha|&\le C G[|Dv|^2]|\nabla D v|^2 &&  \mbox{if} \quad |\alpha| =2,\\
|E_\alpha|&\le C G[|Dv|^2]( |\nabla Dv|^3 + |\nabla^2 Dv| |\nabla Dv|) &&\mbox{if} \quad |\alpha|=3.
\end{aligned}
\eq
\end{lemma}

\begin{lemma}[Coercivity estimates]\label{lem:kkk2}
For any function $G : \R_+ \to \R_+$ satisfying \eqref{H1}--\eqref{H3}, the following hold:

\begin{enumerate}
\item
For any $3\times 3$ matrices $A$ and $B$,  
\begin{equation}\label{eq:kkk2}
G[|A|^2]|B|^2 + 2G'[|A|^2] (A:B)^2 \ge m_0 |B|^2,
\end{equation}
where $m_0$ is from \eqref{H2}.

\item
For any $w_1, w_2 \in H^1(\T^3)$,

\[
m_0 \|D(w_1 -w_2)\|_{L^2}^2 \le  \intt (G[|Dw_1|^2]Dw_1 - G[|Dw_2|^2]Dw_2): D(w_1 -w_2)\,dx,
\]
where $m_0$ is from \eqref{H2}.
\end{enumerate}
\end{lemma}

Finally, we present the local-in-time existence and uniqueness of strong solutions to \eqref{main_sys}-\eqref{init_data}. Since our work is mainly concerned with the global-in-time existence theory, we postpone its proof to Appendix \ref{app.A}.

\begin{theorem}\label{thm:local_ex}
Suppose that the initial data $(\rho_0, u_0, v_0)$ satisfy
\begin{enumerate}
\item $\displaystyle \inf_{x \in \T^3}\rho_0(x) >0$  and
\item $(\rho_0, u_0, v_0) \in H^2(\T^3)\times H^3(\T^3)\times H^3(\T^3)$ with $\nabla\cdot v_0=0$.
\end{enumerate}

Then for any $0<\e_0 < \delta_0 $, there exists a constant $T^*>0$ depending only on $\e_0$ and $\delta_0$ such that if 
\[
\|\rho_0 -\rho_c(0)\|_{H^2} + \|u_0 - u_c(0)\|_{H^3} + \|v_0 - v_c(0)\|_{H^3} + |u_c(0)  - v_c(0)|< \e_0
\]
then system \eqref{main_sys}-\eqref{init_data} admits a unique strong solution 
\[
(\rho, u, v) \in \calC([0,T^*]; H^2(\T^3)) \times \calC([0,T^*]; H^3(\T^3)) \times \calC([0,T^*]; H^3(\T^3))
\]
satisfying $\displaystyle \inf_{(t,x) \in [0,T^*]\times \T^3}\rho(t,x) >0$ and
\[
\sup_{t \in [0,T^*]} \lt( \|\rho(t, \cdot) -\rho_c(t)\|_{H^2} + \|u(t,\cdot) - u_c(t)\|_{H^3} + \|v(t,\cdot) - v_c(t)\|_{H^3} + |u_c(t) - v_c(t)|\rt)<\delta_0.
\]

\end{theorem}

\section{A priori large-time behavior estimates}
\setcounter{equation}{0}

In this section, we present a priori estimates for the large-time behavior of regular solutions to \eqref{main_sys}-\eqref{init_data}. First, we present estimates of the mass, total momentum and energy for the system. For the notation, we recall \eqref{aver_quant} and additionally, we set
\[
m_c(t) := \intt \rho u\,dx.
\]

\begin{proposition}
Let $(\rho, u,v)$ be a solution to \eqref{main_sys}-\eqref{init_data} with sufficient regularity. Then we have the following:

\begin{enumerate}
\item
The mass is preserved in time:

\[
\rho_c(t) = \rho_c(0) = \intt \rho_0 (x)\,dx, \quad \forall t\ge 0.
\]

\item
The total momentum is preserved in time:

\[
m_c(t) + v_c(t) = m_c(0) + v_c(0), \quad \forall t \ge0.
\]
\item
The total energy monotonically decreases in time:
\[
\frac12\frac{d}{dt}\lt(\intt \rho |u|^2\,dx + \intt |v|^2\,dx\rt) + \intt G[|Dv|^2] |Dv|^2\,dx + \intt \rho|u-v|^2\,dx  = 0.
\]
\end{enumerate}
\end{proposition}

\begin{proof}
We can deduce the conservation of the mass and total momentum directly from the system. For the total energy, straightforward computation gives
 
\bq\label{eq:npe_en}
\frac12\frac{d}{dt} \intt \rho |u|^2\,dx = \intt \rho u \cdot (v-u)\,dx,
\eq
and

\bq\label{eq:nn_en}
\frac12\frac{d}{dt}\intt |v|^2\,dx + \intt G[|Dv|^2] |Dv|^2\,dx  = - \intt \rho v \cdot (v-u)\,dx.
\eq
Then we combine \eqref{eq:npe_en} with \eqref{eq:nn_en} to yield the desired estimates for the total energy.
\end{proof}

By conservation of mass,
\bq\label{mass_cons}
\rho_c(t)=\rho_c(0)=1, \quad \forall t \ge 0.
\eq

Since the system is given in the periodic domain, we can expect that the velocity fields would converge toward their averaged quantities. So, we introduce the following modulated energy functional as

\begin{equation}
\calE(t):=\tfrac12\intt\rho\,\abs{u-m_c}^2\,dx+\tfrac12\intt\abs{v-v_c}^2\,dx+\tfrac14\abs{m_c-v_c}^2 .
\end{equation}
Then, the following relation is obvious:

\begin{lemma}\label{lem:l2en}
For $T>0$, let $(\rho, u,v)$ be a solution to \eqref{main_sys}-\eqref{init_data}  on $[0,T]$ with sufficient regularity. Then we have

\[
\frac{d}{dt}\calE(t) + \calD(t) = 0,\quad \forall t>0, 
\] 
where the dissipation functional $\calD = \calD(t)$ is given by
\[
\calD(t):=\intt\Gv\abs{\Dsym v}^2\,dx+\intt\rho\abs{u-v}^2\,dx.
\]
\end{lemma}

To show the temporal decay of the modulated energy functional $\calE$, one needs to show that the dissipation $\calD$ bounds $\calE$. The proof is reminiscent of \cite[Lemma 3.6]{CK16}, but we present it here for readers' convenience.

\begin{lemma}\label{lem:dissbd}
For $T>0$, let $(\rho, u,v)$ be a solution to \eqref{main_sys}-\eqref{init_data} on $[0,T]$ with sufficient regularity. Moreover, assume that there exists a positive constant $\bar\rho>0$ independent of $T$ such that
\[
0 \le \rho(t,x) \le \bar\rho, \quad \forall (t,x) \in [0,T]\times \T^3.
\]
Then, there is a constant $C^*>0$, independent of $T$ such that
\[
\calE \le C^* \calD, \quad \forall t\in[0,T].
\]
\end{lemma}

\begin{proof}
First, by using Poincar\'e inequality,  Korn's inequality  and \eqref{H1}, we obtain
\[ c_0\intt |v-v_c|^2\,dx\le \intt\Gv\abs{\Dsym v}^2\,dx\le \calD,  \]
for some $c_0>0$.

Due to \eqref{mass_cons}, we have
 \[
\intt\rho(u-m_c)\,dx=m_c-m_c\intt\rho\,dx=0.
\]
So the decomposition
\[
u-m_c=(u-v)+(v-v_c)+(v_c-m_c)
\]
provides
\[
m_c-v_c=\intt\rho(u-v)\,dx+\intt\rho(v-v_c)\,dx.
\]
Thus, by Young's inequality and Cauchy-Schwarz inequality,
\[|m_c-v_c|^2\le2\intt\rho\abs{u-v}^2+2\intt\rho\abs{v-v_c}^2\le2\calD+2\bar\rho\intt\abs{v-v_c}^2\le C\calD, \]
for some positive constant $C$ independent of $t$. Similarly, we can get
\[ \intt\rho\abs{u-m_c}^2\le3\intt\rho\abs{u-v}^2+3\intt\rho\abs{v-v_c}^2+3\abs{m_c-v_c}^2\le C\calD, \]
and this implies the desired result.
\end{proof}
Therefore, we combine Lemma \ref{lem:l2en} with \ref{lem:dissbd} to deduce the following large-time behavior estimate 

\begin{proposition}\label{prop:decay}
For $T>0$, let $(\rho, u,v)$ be a solution to \eqref{main_sys}-\eqref{init_data} on $[0,T]$ with sufficient regularity. Moreover, assume that there exists a positive constant $\bar\rho>0$ independent of $T$ such that
\[
0 \le \rho(t,x) \le \bar\rho, \quad \forall (t,x) \in [0,T]\times \T^3.
\]
Then there exists a positive constant $\theta>0$ independent of $T$ such that
\[
\calE(t) \le \calE(0)e^{-\theta t}, \quad \forall t\in[0,T].
\]
\end{proposition}

\begin{remark}\label{rem:mc_uc}
Note that the mean velocity $u_c(t)$ and mean momentum $m_c(t)$ do not necessarily coincide. However, under the additional assumption that $\rho(t,x) \ge \underline{\rho}>0$, their difference can be controlled by the modulated energy. In particular, for global-in-time smooth solutions, they become asymptotically identical as $t \to \infty$. Indeed, using  $\rho_c(t) \equiv 1$, we have
\[\begin{aligned}
|m_c(t) -  u_c(t)| &= \lt|\int_{\T^3} (\rho-1)(y)u(y)\,dy\rt|\\
&= \lt|\int_{\T^3} (\rho-1)(y)(u(y)- m_c(t))\,dy\rt|\\
&\le (1+\bar\rho)\|u-m_c\|_{L^2}\le \frac{1+\bar\rho}{\sqrt{\underline{\rho}}}\sqrt{\mathcal{E}(t)}.
\end{aligned}\]
On the other hand, replacing $m_c(t)$ by $u_c(t)$ in the second equality yields
\[
|m_c(t)-u_c(t)| \le (1+\bar\rho)\|u-u_c(t)\|_{L^2},
\]
which also provides a direct estimate in terms of the velocity fluctuation.
\end{remark}

\section{A priori energy estimates}

In this section, we study a priori estimates for $L^2$-Sobolev norms of regular solutions to \eqref{main_sys}-\eqref{init_data}. Since we have assumed $\rho_c (t) \equiv 1$, we introduce a new variable $h := \rho-1$ to rewrite \eqref{main_sys} as follows:

\[
\begin{aligned}
&\pa_t h +\nabla \cdot ((1+h) u)=0, \\
&\pa_t u+(u \cdot \nabla) u =(v-u),\\
&\pa_t v +(v\cdot\nabla)v+\nabla p - \nabla\cdot \!\bigl(\Gv\,\Dsym v\bigr)=-(1+h)(v-u),\\
&\nabla\cdot v=0,
\end{aligned}
\]

 One notes that the modulated energy estimates in Lemma \ref{lem:l2en} already provide  $L^2$-estimates of solutions. Thus, it suffices to investigate the $L^2$-norms of their derivatives. We first consider the incompressible fluid.

\begin{lemma}\label{lem:l2est_v}
Let $(\rho, u,v)$ be a solution to \eqref{main_sys}-\eqref{init_data} on the time interval $[0,T]$ with sufficient regularity. Then there exist positive constants $c_0$ and $C$ independent of $t$ satisfying
\bq\label{eq:1st_l2v}
\begin{aligned}
\frac{d}{dt}\|\nabla v\|_{L^2}^2 + 2c_0 \|\nabla^2 v\|_{L^2}^2 &\le C\|\nabla v\|_{L^\infty} \|\nabla v\|_{L^2}^2 \\
&\quad+ C(1+\|h\|_{L^\infty})^{1/2} \|(1+h)^{1/2} (u-v)\|_{L^2}\|\nabla^2 v\|_{L^2},
\end{aligned}
\eq
and for $\ell=2,3$,

\bq\label{eq:high_l2v}
\begin{aligned}
\frac{d}{dt}&\|\nabla^\ell v\|_{L^2}^2 + 2c_0 \|\nabla^{\ell+1} v\|_{L^2}^2 \\
&\le C\|\nabla v\|_{L^\infty}\|\nabla^\ell v\|_{L^2}^2\\
&\quad + C\left(\max_{0\le s\le C\|v-v_c\|_{H^3}^2}G[s]\right)(1+\|v-v_c\|_{H^3}) \|v-v_c\|_{H^3}^{3/2} \|\nabla^{\ell+1} v\|_{L^2}^{3/2}\\
&\quad +  C\lt(\|\nabla^{\ell-1} h\|_{L^2} \|u-v\|_{L^\infty} + (1+\|h\|_{L^\infty})\|\nabla^{\ell-1}(u-v)\|_{L^2}\rt)\|\nabla^{\ell+1}v\|_{L^2}.
\end{aligned}
\eq
\end{lemma}

\begin{proof}
We split the proof into two parts:\\

\noindent $\bullet$ (Part A: first-order estimates) Direct computation yields, for $i=1,2,3$,

\[
\begin{aligned}
\frac12\frac{d}{dt}\|\pa_i v\|_{L^2}^2 &= \intt \pa_i \nabla\cdot (G[|Dv|^2] Dv) \cdot \pa_i v\,dx\\
&\quad - \intt \pa_i [(v \cdot \nabla) v] \cdot \pa_i v\,dx + \intt (1+h)(v-u)\cdot\pa_i^2 v\,dx\\
&=: \sfI_{11} + \sfI_{12} + \sfI_{13}.
\end{aligned}
\]
For $\sfI_{11}$, note that
\[ \pa_i(\Gv\Dsym v)=\Gv\,\pa_i\Dsym v+2G'[\abs{\Dsym v}^2]\,(\Dsym v:\pa_i\Dsym v)\,\Dsym v . \]
Hence, we use \eqref{eq:kkk2} to get
\[\begin{aligned}
\sfI_{11} &= -\intt \pa_i (G[|Dv|^2]Dv) : \pa_i Dv\,dx \le -m_0 \|\pa_i Dv\|_{L^2}^2 
\end{aligned}\]
For $\sfI_{12}$, 
\[\begin{aligned}
\sfI_{12} &= -\intt [(\pa_i v \cdot \nabla) v] \cdot \pa_i v\,dx\\
&\le C\|\nabla v\|_{L^\infty}\|\pa_i v\|_{L^2}^2.
\end{aligned}\]
For $\sfI_{13}$, we have
\[\begin{aligned}
\sfI_{13} \le (1+\|h\|_{L^\infty})^{1/2} \| (1+h)^{1/2} (u-v)\|_{L^2}\|\pa_i^2 v\|_{L^2}.
\end{aligned}\]
Then we collect the estimates for $\sfI_{1j}$'s, sum over $i=1,2,3$ and use Korn's inequality to get

\[
\begin{aligned}
\frac{d}{dt}\|\nabla v\|_{L^2}^2 + 2c_1 \|\nabla^2 v\|_{L^2}^2 &\le C\|\nabla v\|_{L^\infty} \|\nabla v\|_{L^2}^2 \\
&\quad+ C(1+\|h\|_{L^\infty})^{1/2} \|(1+h)^{1/2} (u-v)\|_{L^2}\|\nabla^2 v\|_{L^2},
\end{aligned}
\]
for some positive constants $c_1$ and $C$ independent of $T$.\\

\noindent $\bullet$ (Part B: second- and third-order estimates) For a multi-index $\beta$ with $|\beta|=1$ or 2, and $i=1,2,3$, one has

\[\begin{aligned}
\frac12\frac{d}{dt}\|\pa^\beta \pa_i v\|_{L^2}^2 &= \intt \pa^\beta \pa_i \nabla\cdot (G[|Dv|^2]Dv)\cdot \pa^\beta \pa_i v\,dx\\
&\quad - \intt \pa^\beta \pa_i [(v \cdot \nabla) v] \cdot \pa^\beta \pa_i v\,dx + \intt \pa^\beta [(1+h)(v-u)]\cdot\pa^\beta \pa_i^2 v\,dx\\
&=: \sfI_{21} + \sfI_{22} + \sfI_{23}.
\end{aligned}\]
For $\sfI_{21}$, one uses Lemma \ref{lem:kkk1}, Sobolev embedding and Korn's inequality to get if $|\beta|=1$,

\[\begin{aligned}
\sfI_{21} &= -\intt \pa^\beta \pa_i (G[|Dv|^2]Dv) :\pa^\beta \pa_i Dv\,dx\\
&\le -m_0 \|\pa^\beta \pa_i Dv\|_{L^2}^2 + C\left(\max_{0\le s\le \|\nabla v\|_{L^\infty}^2}G[s]\right)\|\nabla Dv\|_{L^6}\|\nabla Dv\|_{L^3}\|\pa^\beta \pa_i Dv\|_{L^2}\\
&\le -m_0 \|\pa^\beta \pa_i Dv\|_{L^2}^2 + C\left(\max_{0\le s\le C\|v-v_c\|_{H^3}^2}G[s]\right)\|v-v_c\|_{H^3}^2 \|\nabla^3 v\|_{L^2},
\end{aligned}\]
and if $|\beta|=2$,

\[\begin{aligned}
\sfI_{21} &\le  -m_0 \|\pa^\beta \pa_i Dv\|_{L^2}^2\\
&\quad + C\left(\max_{0\le s\le \|\nabla v\|_{L^\infty}^2}G[s]\right)\lt( \|\nabla Dv\|_{L^6}^3 + \|\nabla^2 Dv\|_{L^3}\|\nabla Dv\|_{L^6}\rt)\|\pa^\beta \pa_i Dv\|_{L^2}\\
&\le  -m_0 \|\pa^\beta \pa_i Dv\|_{L^2}^2 + C\left(\max_{0\le s\le C\|v-v_c\|_{H^3}^2}G[s]\right)(1+\|v-v_c\|_{H^3})\|v-v_c\|_{H^3}^{3/2} \|\nabla^4 v\|_{L^2}^{3/2},
\end{aligned}\]
where we used $\|\nabla^3 f \|_{L^3} \le C\|\nabla^3 f\|_{L^2}^{1/2} \|\nabla^4 f\|_{L^2}^{1/2}$.

In either case, we obtain

\[
\sfI_{21} \le  -m_0 \|\pa^\beta \pa_i Dv\|_{L^2}^2 + C\left(\max_{0\le s\le C\|v-v_c\|_{H^3}^2}G[s]\right)(1+\|v-v_c\|_{H^3})\|v-v_c\|_{H^3}^{3/2} \|\nabla^{|\beta|+2} v\|_{L^2}^{3/2}.
\]
For $\sfI_{22}$ and $\sfI_{23}$, one uses Moser-type estimates to yield

\[\begin{aligned}
\sfI_{22} &\le C\|\nabla v\|_{L^\infty}\|\nabla^{|\beta|+1}v\|_{L^2}^2,
\end{aligned}\]
and
\[\begin{aligned}
\sfI_{23} &\le C\lt(\|\nabla^{|\beta|} h\|_{L^2} \|u-v\|_{L^\infty} + (1+\|h\|_{L^\infty})\|\nabla^{|\beta|}(u-v)\|_{L^2}\rt)\|\nabla^{|\beta|+2}v\|_{L^2}.
\end{aligned}\]
Hence, by gathering the estimates for $\sfI_{2j}$'s and summing over all $i=1,2,3$ and $\beta$, we can find a positive constant $c_2>0$ independent of $T$ such that

\[
\begin{aligned}
\frac{d}{dt}&\|\nabla^\ell v\|_{L^2}^2 + 2c_2 \|\nabla^{\ell+1} v\|_{L^2}^2 \\
&\le C\|\nabla v\|_{L^\infty}\|\nabla^\ell v\|_{L^2}^2\\
&\quad + C\left(\max_{0\le s\le C\|v-v_c\|_{H^3}^2}G[s]\right)(1+\|v-v_c\|_{H^3}) \|v-v_c\|_{H^3}^{3/2} \|\nabla^{\ell+1} v\|_{L^2}^{3/2}\\
&\quad +  C\lt(\|\nabla^{\ell-1} h\|_{L^2} \|u-v\|_{L^\infty} + (1+\|h\|_{L^\infty})\|\nabla^{\ell-1}(u-v)\|_{L^2}\rt)\|\nabla^{\ell+1}v\|_{L^2},
\end{aligned}
\]
for $\ell=2,3$ and choosing $c_0 := \min\{c_1, c_2\}$  implies the desired estimates.
\end{proof}

Next, we proceed to the estimates for the velocity field $u$ of the pressureless flow.

\begin{lemma}\label{lem:l2est_u}
Let $(\rho, u,v)$ be a solution to \eqref{main_sys}-\eqref{init_data} on the time interval $[0,T]$ with sufficient regularity. Then there exists a positive constant $C$ independent of $t$ satisfying
\[
\frac{d}{dt}\|\nabla^\ell u\|_{L^2}^2 + \|\nabla^\ell u\|_{L^2}^2 \le C\|\nabla u\|_{L^\infty} \|\nabla^\ell u\|_{L^2}^2 + \|\nabla^\ell v\|_{L^2}^2,
\]
for each $1 \le \ell \le 3$.
\end{lemma}

\begin{proof}
For a multi-index $\alpha$ with $1\le |\alpha|\le 3$, by Moser-type estimates and Young's inequality, 

\[\begin{aligned}
\frac12\frac{d}{dt}\|\pa^\alpha u\|_{L^2}^2 + \|\pa^\alpha u\|_{L^2}^2 &= - \intt \pa^\alpha [(u \cdot \nabla) u] \cdot \pa^\alpha u\,dx\\
&\quad +\intt \pa^\alpha v \cdot \pa^\alpha u\,dx\\
&=  \frac12 \intt (\nabla \cdot u) |\pa^\alpha u|^2\,dx - \intt \lt(\pa^\alpha [(u \cdot \nabla )u] - (u \cdot \nabla) \pa^\alpha u \rt) \cdot \pa^\alpha u\,dx\\
&\quad + \intt \pa^\alpha v \cdot \pa^\alpha u\,dx\\
&\le C\|\nabla u\|_{L^\infty}\|\nabla^{|\alpha|} u\|_{L^2}^2 + \frac12\|\pa^\alpha v\|_{L^2}^2  + \frac12\|\pa^\alpha u\|_{L^2}^2,
\end{aligned}\]
and summing the above relation over $\alpha$ implies the desired result.

\end{proof}

Finally, we present the estimates for the density of the pressureless flow.

\begin{lemma}\label{lem:l2est_h}
Let $(\rho, u,v)$ be a solution to \eqref{main_sys}-\eqref{init_data} on the time interval $[0,T]$ with sufficient regularity. Then there exists a positive constant $C$ independent of $t$ satisfying

\bq\label{eq:h2_h}
\begin{aligned}
\frac{d}{dt}\|h\|_{H^2} \le C\lt(\|\nabla u\|_{L^\infty} + \|\nabla u\|_{H^2}\rt)\|h\|_{H^2} + C\|\nabla u\|_{H^2}.
\end{aligned}
\eq
\end{lemma}

\begin{proof}
For a multiindex $\beta$ with $0 \le |\beta| \le 2$, by Moser-type estimates,

\[\begin{aligned}
\frac12\frac{d}{dt}\|\pa^\beta h\|_{L^2}^2 &= -\intt \pa^\beta \nabla\cdot (hu) \pa^\beta h\,dx - \intt \pa^\beta h \nabla\cdot \pa^\beta u\,dx\\
&= \frac12\intt (\nabla\cdot u) |\pa^\beta h|^2\,dx - \intt \lt(\pa^\beta \nabla\cdot (hu) - u \cdot \pa^\beta\nabla h \rt) \pa^\beta h\,dx\\
&\quad -\intt \pa^\beta h \pa^\beta \nabla\cdot u\,dx\\
&\le C\|\nabla\cdot u\|_{L^\infty}\|\pa^\beta h\|_{L^2}^2 + C\|\pa^\beta h\|_{L^2}\lt(\|h\|_{L^\infty}\|\nabla^{|\beta|+1} u\|_{L^2} + \|\nabla u\|_{L^\infty}\|\nabla^{|\beta|} h\|_{L^2} \rt)\\
&\quad + \|\pa^\beta h\|_{L^2} \|\pa^\beta \nabla\cdot u\|_{L^2}.
\end{aligned}\]
Summing over all $\beta$, dividing by $(\|h\|_{H^2}^2+\delta^2)^{1/2}$, and letting $\delta\downarrow0$,

\[\begin{aligned}
\frac{d}{dt}\|h\|_{H^2} \le C\lt(\|\nabla u\|_{L^\infty} + \|\nabla u\|_{H^2}\rt)\|h\|_{H^2} + C\|\nabla u\|_{H^2},
\end{aligned}\]
which is the desired estimate.

\end{proof}

\section{Proof of Theorem \ref{main_thm}}
\setcounter{equation}{0}
Now, we turn to the proof of Theorem \ref{main_thm}. 

\subsection{Global-in-time existence}

We proceed by a bootstrapping argument, based on a priori large-time behavior and energy estimates. First, we set
\[
\mathcal{N}(t) := \|h(t)\|_{H^2} + \|(u -u_c)(t)\|_{H^3} + \|(v-v_c)(t)\|_{H^3} + |u_c(t)-v_c(t)|,
\]
and assume that
\bq\label{eq:boots}
\sup_{0\le t \le T} \mathcal{N}(t) < \eta
\eq
for sufficiently small constant $\eta>0$. Then we establish the uniform estimates for velocity fields.

\begin{proposition}\label{prop:unif_t_est}
For $T>0$, let $(\rho, u, v)$ be a solution to \eqref{main_sys}-\eqref{init_data} on the time interval $[0,T]$ with sufficient regularity and satisfying the bootstrap bound \eqref{eq:boots} with sufficiently small $\eta>0$. Then there exists a positive constant $C_1^*>1$ which is independent of $T$ and $\mathcal{N}(0)$ satisfying
\[
\sup_{0 \le t \le T}\lt(\|(u -u_c)(t)\|_{H^3} + \|(v-v_c)(t)\|_{H^3}+ |u_c(t)-v_c(t)|\rt) \le C_1^* \mathcal{N}(0).
\]
\end{proposition}
\begin{proof}
First, we recall from Lemma \ref{lem:l2en} that
\bq\label{eq:en_ineq}
\mathcal{E}(t)  +\int_0^t \mathcal{D}(\tau)\,d\tau \le \mathcal{E}(0), \quad \forall t \in [0,T].
\eq

Since the bootstrap bound \eqref{eq:boots} gives $\frac12\le\rho\le\frac32$, we can use the arguments in Remark \ref{rem:mc_uc} to have
\bq\label{eq:en_bd_new}
\begin{aligned}
\mathcal{E}(0)&\le C\lt(\|u_0-u_c(0)\|_{L^2}^2+\|v_0-v_c(0)\|_{L^2}^2+|u_c(0)-v_c(0)|^2\rt)\\
&\le C(\mathcal{N}(0))^2.
\end{aligned}
\eq
Thus, \eqref{eq:en_ineq} implies

\bq\label{eq:diss_bd1}
\int_0^T \|\nabla v(t)\|_{L^2}^2\,dt + \int_0^T \intt (1+h)|u-v|^2\,dxdt \le C(\mathcal{N}(0))^2,
\eq
for some constant $C>0$ independent of $T$. By applying Poincar\'e inequality, Young's inequality and the bootstrap bound \eqref{eq:boots} to \eqref{eq:1st_l2v}, we get

\[\begin{aligned}
\frac{d}{dt}\|\nabla v\|_{L^2}^2 + 2c_0\|\nabla^2 v\|_{L^2}^2 &\le C\eta\|\nabla^2 v\|_{L^2}^2 +  C\lt( \intt (1+h) |u-v|^2\,dx\rt)^{1/2}\|\nabla^2 v\|_{L^2}\\
&\le c_0 \|\nabla^2 v\|_{L^2}^2 + C\intt (1+h) |u-v|^2\,dx,
\end{aligned}\]
and hence, we use \eqref{eq:diss_bd1} to get

\bq\label{eq:h1_bd_v}
\begin{aligned}
\sup_{0 \le t \le T}\|\nabla v(t)\|_{L^2}^2 + \int_0^T \|\nabla^2 v (t)\|_{L^2}^2\,dt &\le C\|\nabla v_0\|_{L^2}^2 + C\int_0^T \intt (1+h) |u-v|^2\,dxdt\\
&\le C(\mathcal{N}(0))^2.
\end{aligned}
\eq
Next, we apply the bootstrap bound \eqref{eq:boots} to Lemma \ref{lem:l2est_u} with $\ell=1,2$ to get

\[\begin{aligned}
\frac{d}{dt}\|\nabla u\|_{H^1}^2 + \frac12\|\nabla u\|_{H^1}^2 \le C\|\nabla v\|_{H^1}^2,
\end{aligned}\]
and one combines \eqref{eq:diss_bd1} with \eqref{eq:h1_bd_v} to yield

\bq\label{eq:h2_bd_u}
\begin{aligned}
\sup_{0\le t \le T} \|\nabla u(t)\|_{H^1}^2 + \int_0^T \|\nabla u(t)\|_{H^1}^2\,dt &\le C\|\nabla u_0\|_{H^1}^2 + C\int_0^T \|\nabla v(t)\|_{H^1}^2\,dt\\
&\le C(\mathcal{N}(0))^2.
\end{aligned}
\eq

Now, by applying the bootstrap bound, Poincare inequality and Sobolev embedding $H^2 \hookrightarrow L^\infty$ to \eqref{eq:high_l2v}

\[\begin{aligned}
&\frac{d}{dt}\|\nabla^2 v\|_{H^1}^2 + 2c_0 \|\nabla^3 v\|_{H^1}^2 \\
&\le C\eta\|\nabla^3 v\|_{H^1}^2 + C\lt(\|u-v\|_{L^\infty} + \|\nabla u\|_{H^1} + \|\nabla v\|_{H^1}\rt)\|\nabla^3 v\|_{H^1}\\
&\le C\eta\|\nabla^3 v\|_{H^1}^2+ C\lt(\|u-v\|_{L^2} + \|\nabla u\|_{H^1} + \|\nabla v\|_{H^1}\rt)\|\nabla^3 v\|_{H^1}.
\end{aligned}\]
Here, note that \eqref{eq:boots} implies
\[
\frac12 \le (1+h)(t,x) \le \frac32, \quad \forall (t, x) \in [0,T]\times\T^3
\]
and in turn,
\[
\|u-v\|_{L^2} \le C\|(1+h)^{1/2}(u-v)\|_{L^2}.
\]
Thus, we gather \eqref{eq:diss_bd1}, \eqref{eq:h1_bd_v} and \eqref{eq:h2_bd_u} and use Young's inequality to have

\bq\label{eq:h3_bd_v}
\begin{aligned}
\sup_{0\le t \le T} \|\nabla^2 v(t)\|_{H^1}^2 + \int_0^T \|\nabla^3 v(t)\|_{H^1}^2\,dt &\le C\|\nabla^2 v_0\|_{H^1}^2 + C\int_0^T \intt (1+h)|u-v|^2\,dxdt \\
&\quad + C\int_0^T \lt(\|\nabla v(t)\|_{H^1}^2 + \|\nabla u(t)\|_{H^1}^2\rt)\,dt\\
&\le C(\mathcal{N}(0))^2.
\end{aligned}
\eq
Similarly to \eqref{eq:h2_bd_u}, we can also get

\bq\label{eq:h4_bd_u}
\begin{aligned}
\sup_{0\le t \le T}\|\nabla^3 u(t)\|_{L^2}^2 + \int_0^T \|\nabla^3 u(t)\|_{L^2}^2 \,dt &\le C\|\nabla^3 u_0\|_{L^2}^2 + C\int_0^T \|\nabla^3 v(t)\|_{L^2}^2\,dt\\
&\le C(\mathcal{N}(0))^2.
\end{aligned}
\eq
Finally, Proposition \ref{prop:decay} and Remark \ref{rem:mc_uc} imply
\bq\label{eq:mean_diff}
|u_c(t)-v_c(t)|\le |m_c(t)-v_c(t)|+|m_c(t)-u_c(t)|
\le C\sqrt{\mathcal{E}(t)}\le C\mathcal{N}(0)e^{-\theta t/2}.
\eq
Hence, by collecting \eqref{eq:en_ineq}, \eqref{eq:h1_bd_v}, \eqref{eq:h2_bd_u}, \eqref{eq:h3_bd_v}, \eqref{eq:h4_bd_u}, and \eqref{eq:mean_diff}, one has

\bq\label{eq:unif_uv}
\sup_{0\le t \le T} \lt(\|(u-u_c)(t)\|_{H^3} + \|(v-v_c)(t)\|_{H^3}+ |u_c(t)-v_c(t)|\rt) \le C_1^*\mathcal{N}(0),
\eq
for some positive constant $C_1^*>1$ independent of $T$ and this implies the desired result.

\end{proof}

Based on the previous uniform estimates, we show the exponential decay rate of velocity fluctuations in Sobolev spaces.

\begin{corollary}\label{cor:LT_high}
For $T>0$, let $(\rho, u, v)$ be a solution to \eqref{main_sys}-\eqref{init_data} on the time interval $[0,T]$ with sufficient regularity and satisfying the bootstrap bound \eqref{eq:boots} with sufficiently small $\eta>0$. Then there exist positive constants $C, c>0$ which are independent of $T$ satisfying
\[
\lt(\|(u -u_c)(t)\|_{H^3} + \|(v-v_c)(t)\|_{H^3}\rt) \le C \mathcal{N}(0) e^{-c t}, \quad \forall t \in [0,T].
\]
\end{corollary}

\begin{proof}
First, Proposition \ref{prop:decay}, together with Remark \ref{rem:mc_uc} and \eqref{eq:en_bd_new}, implies the exponential decay rate for $L^2$-norms of $u-u_c$ and $v-v_c$. For exponential decay rate of derivatives, we use Proposition \ref{prop:decay} and \eqref{eq:unif_uv} to obtain a positive constant $\theta_1$ independent of $T$ satisfying
\bq\label{eq:u_decay}
\begin{aligned}
\|\nabla u\|_{L^\infty} &\le C\|\nabla^3 u\|_{L^2}^{5/6}\|u-u_c\|_{L^2}^{1/6}\le C\mathcal{N}(0) e^{-\theta_1 t},\\
\|\nabla^\ell u\|_{L^2} &\le C\|\nabla^3 u\|_{L^2}^{\frac{\ell}{3}} \|u-u_c\|_{L^2}^{\frac{3-\ell}{3}} \le C\mathcal{N}(0) e^{-\theta_1 t},
\end{aligned}
\eq
and similarly,
\bq\label{eq:v_decay}
\begin{aligned}
\|\nabla v\|_{L^\infty} &\le  C\mathcal{N}(0) e^{-\theta_1 t},\\
\|\nabla^\ell v\|_{L^2} &\le  C\mathcal{N}(0) e^{-\theta_1 t},
\end{aligned}
\eq
for $\ell = 1,2$.

Then we apply Poincar\'e inequality to \eqref{eq:high_l2v} and obtain

\[\begin{aligned}
&\frac{d}{dt} \|\nabla^3 v\|_{L^2}^2 + 2c_0\|\nabla^4 v\|_{L^2}^2 \\
&\le C\eta \|\nabla^4 v\|_{L^2}^2 + C\lt( \|u-v\|_{L^\infty} + \|\nabla^2 (u-v)\|_{L^2}\rt)\|\nabla^4 v\|_{L^2}\\
&\le C\eta \|\nabla^4 v\|_{L^2}^2 + C\lt( \|u-v\|_{L^2} + \|\nabla^2 u\|_{L^2} + \|\nabla^2 v\|_{L^2}\rt)\|\nabla^4 v\|_{L^2}.
\end{aligned}\]
Since
\[\begin{aligned}
\|u-v\|_{L^2}^2 &\le C\intt \rho |u-v|^2\,dx = C\intt \rho |u-m_c+m_c - v_c + v_c -v|^2\,dx  \\
&\le C\lt(\intt \rho|u-m_c|^2 + |m_c - v_c|^2 + \intt \rho| v- v_c|^2\,dx\rt)\\
&\le C\mathcal{E}(t),
\end{aligned}\]
we gather Proposition \ref{prop:decay}, \eqref{eq:u_decay} and \eqref{eq:v_decay} to yield

\[\begin{aligned}
&\frac{d}{dt} \|\nabla^3 v\|_{L^2}^2 + 2c_0\|\nabla^4 v\|_{L^2}^2 \\
&\le C\eta \|\nabla^4 v\|_{L^2}^2 + C\mathcal{N}(0)e^{-\theta_2 t}\|\nabla^4 v\|_{L^2}\\
&\le C\eta \|\nabla^4 v\|_{L^2}^2 + \frac{c_0}{2}\|\nabla^4 v\|_{L^2}^2 +  C\lt(\mathcal{N}(0)\rt)^2e^{-2\theta_2 t},
\end{aligned}\]
for some positive constant $\theta_2$ independent of $T$, where we used Young's inequality. Hence, the smallness of $\eta$ and the Poincar\'e inequality imply

\[
\frac{d}{dt}\|\nabla^3 v\|_{L^2}^2 + \tilde{c}_1 \|\nabla^3 v\|_{L^2}^2 \le C\lt(\mathcal{N}(0)\rt)^2e^{-2\theta_2 t},
\]
for some positive constant $\tilde{c}_1>0$ and Gr\"onwall's lemma gives
\bq\label{eq:v3_exp_decay}
\|\nabla^3 v(t)\|_{L^2}^2 \le C\lt(\mathcal{N}(0)\rt)^2e^{-2\theta_3 t},
\eq
for any $t\in [0,T]$, where $\theta_3>0$ is a positive constant independent of $T$.  Finally, we combine Lemma \ref{lem:l2est_u} with \eqref{eq:v3_exp_decay} to attain

\[\begin{aligned}
\frac{d}{dt}\|\nabla^3 u\|_{L^2}^2 + \|\nabla^3 u\|_{L^2}^2 &\le C\eta \|\nabla^3 u\|_{L^2}^2 + \|\nabla^3 v\|_{L^2}^2\\
&\le \frac12\|\nabla^3 u\|_{L^2}^2 + C\lt(\mathcal{N}(0)\rt)^2e^{-2\theta_3 t}
\end{aligned}\]
due to the smallness of $\eta$. Thus,  we use Gr\"onwall's lemma again to obtain

\bq
\|\nabla^3 u(t)\|_{L^2}^2 \le C\lt(\mathcal{N}(0)\rt)^2e^{-2\theta_4 t},
\eq

where $\theta_4 = \frac12\min\{1/4, \theta_3\}$ is a positive constant independent of $T$. Therefore, we gather all the estimates until now to conclude the proof.
\end{proof}

Finally, we present the uniform estimates for the Euler density.

\begin{corollary}\label{cor:h_unif}
For $T>0$, let $(\rho, u, v)$ be a solution to \eqref{main_sys}-\eqref{init_data} on the time interval $[0,T]$ with sufficient regularity and satisfying the bootstrap bound \eqref{eq:boots} with sufficiently small $\eta>0$. Then there exists a positive constant $C_2^*>1$ which is independent of $T$ satisfying
\[
\sup_{0 \le t \le T} \|h(t)\|_{H^2} \le C_2^*\mathcal{N}(0).
\]
\end{corollary}

\begin{proof}

Applying Corollary \ref{cor:LT_high} to \eqref{eq:h2_h}, we obtain

\[\begin{aligned}
\frac{d}{dt}\|h\|_{H^2} \le C\mathcal{N}(0)e^{-c t} \|h\|_{H^2} + C\mathcal{N}(0) e^{-c t}
\end{aligned}\]
and Gr\"onwall's lemma together with $\mathcal{N}(0)<\eta \ll 1$ implies the desired result.
\end{proof}

Now, we are ready to proceed to the proof of our main theorem.

\begin{proof}[Proof of Theorem \ref{main_thm}]
First, we define $C^*>1$ as
\[
C^* =: C_1^* + C_2^*,
\]
where $C_1^*$ and $C_2^*$ are from Proposition \ref{prop:unif_t_est} and Corollary \ref{cor:h_unif}, respectively.

Let $(\rho_0, u_0, v_0)$ be an initial data satisfying
\[
\|\rho_0 -1\|_{H^2} + \|u_0 - u_c(0)\|_{H^3} +\|v_0 -v_c(0)\|_{H^3} + |u_c(0)-v_c(0)| < \e_0,
\] 
where $\e_0$ satisfies $C^*\e_0 <\eta$ for $\eta$ given in \eqref{eq:boots}. Then Theorem \ref{thm:local_ex} guarantees there exists $T^*>0$ such that \eqref{main_sys}-\eqref{init_data} admits  a unique strong solution $(\rho, u, v)$ on the time interval $[0,T^*]$.
satisfying
\[
\sup_{0\le t \le  T^*} \mathcal{N}(t) <\eta.
\]
Then the set $\mathcal{S}$ defined as

\[
\mathcal{S} := \{ T>0 \ : \ \sup_{0\le t \le T} \mathcal{N}(t) < \eta\}
\]
 is nonempty.  If we assume that $\sup \mathcal{S}=:T_{\max} <\infty$, 
we combine Proposition \ref{prop:unif_t_est} with Corollary \ref{cor:h_unif} to yield
\[
\eta = \lim_{t \to T_{\max}}\mathcal{N}(t) \le C^*\e_0 <\eta
\]
which leads to a contradiction and hence, the constructed solution must exist globally in time. The large-time behavior in $H^3(\T^3)$ directly follows from Corollary \ref{cor:LT_high}.
\end{proof}

\section*{Acknowledgment}
 The work of J. Jung was supported by NRF grant no. RS-2026-25589791.

\appendix

\section{Proof of Theorem \ref{thm:local_ex}}\label{app.A}

In this appendix, we prove the local-in-time existence and uniqueness of solutions to \eqref{main_sys}-\eqref{init_data}. We first investigate the uniqueness.

\subsection{Uniqueness} For $i=1,2$ and $T>0$, let $(\rho_i, u_i, v_i)$ be a solution to \eqref{main_sys}-\eqref{init_data} corresponding to the initial data $(\rho_{i,0}, u_{i,0}, v_{i,0}) \in H^2(\T^3) \times H^3(\T^3) \times H^3(\T^3)$ satisfying
\[
(\rho_i, u_i, v_i) \in \calC([0,T];H^2(\T^3))\times \calC([0,T]; H^3(\T^3)) \times \calC([0,T];H^3(\T^3)).
\]
Then we measure $L^2$-differences of these two solutions. We first estimate

\bq\label{eq:udiff_l2}
\begin{aligned}
\frac12\frac{d}{dt}\|u_1 - u_2\|_{L^2}^2 &= -\intt (u_1 - u_2) \cdot [(u_1 - u_2)\cdot \nabla] u_2 \,dx\\
&\quad - \intt (u_1 -u_2)\cdot [u_1 \cdot \nabla(u_1 - u_2)]\,dx  \\
&\quad + \intt (u_1 - u_2) \cdot (v_1 -v_2)\,dx - \|u_1- u_2\|_{L^2}^2\\
&\le C(\|\nabla u_1\|_{L^\infty}+\|\nabla u_2\|_{L^\infty})\|u_1 - u_2\|_{L^2}^2 + \|v_1 - v_2\|_{L^2}^2.
\end{aligned}
\eq
For $i=1,2,3$,  we also obtain,

\bq\label{eq:udiff_h1}
\begin{aligned}
\frac12\frac{d}{dt}\|\pa_i (u_1 - u_2)\|_{L^2}^2 &= - \intt \pa_i (u_1 - u_2)\cdot \pa_i [(u_1 - u_2)\cdot \nabla u_2]\,dx \\
&\quad- \intt \pa_i (u_1 -u_2) \cdot \pa_i [u_1 \cdot \nabla (u_1 -u_2)]\,dx\\
&\quad + \intt \pa_i (u_1 -u_2) \cdot \pa_i (v_1 - v_2)\,dx - \|\pa_i (u_1 - u_2)\|_{L^2}^2\\
&\le  C(\|u_1\|_{H^3}+\|u_2\|_{H^3})\|u_1 - u_2\|_{H^1}^2 + C\|\nabla (v_1 - v_2)\|_{L^2}^2.
\end{aligned}
\eq
Then we first combine \eqref{eq:udiff_l2} with \eqref{eq:udiff_h1} to obtain

\bq\label{eq:udiff_est_comb}
\frac{d}{dt}\|u_1 - u_2\|_{H^1}^2 \le C\lt(\|u_1 - u_2\|_{H^1}^2 + \|v_1 - v_2\|_{L^2}^2\rt) + C_1 \|\nabla (v_1 -v_2)\|_{L^2}^2,
\eq
for some positive constant $C_1>0$. \\

For differences of non-Newtonian fluids, we use Lemma \ref{lem:kkk2} to get

\bq\label{eq:vdiff}
\begin{aligned}
\frac12\frac{d}{dt}\|v_1 - v_2\|_{L^2}^2 &= -\intt (v_1 -v_2)\cdot[ (v_1 - v_2)\cdot \nabla] v_2\,dx\\
&\quad -\intt (G[|Dv_1|^2]D v_1 - G[|Dv_2|^2]Dv_2) : D(v_1- v_2)\,dx\\
&\quad -\intt (v_1 - v_2)\cdot[ (\rho_1 -\rho_2) (v_1 - u_1) + \rho_2 ( v_1 - v_2 - u_1 + u_2)]\,dx\\
&\le C\|\nabla v_2\|_{L^\infty}\|v_1 - v_2\|_{L^2}^2  - m_0 \|D(v_1 - v_2)\|_{L^2}^2\\
&\quad + C(\|v_1\|_{L^\infty} + \|u_1\|_{L^\infty})\|\rho_1 - \rho_2\|_{L^2}\|v_1 -v_2\|_{L^2} \\
&\quad+ C\|\rho_2\|_{L^\infty}\lt(\|v_1- v_2\|_{L^2}^2 + \|v_1 - v_2\|_{L^2}\|u_1 - u_2\|_{L^2}\rt)\\
&\le -\delta \|\nabla (v_1 - v_2)\|_{L^2}^2 + C\lt(\|\rho_1 - \rho_2\|_{L^2}^2 + \|v_1 - v_2\|_{L^2}^2 + \|u_1 - u_2\|_{L^2}^2 \rt),
\end{aligned}
\eq
for some $\delta>0$. 

Finally, we have
\bq\label{eq:rhodiff}\begin{aligned}
\frac12\frac{d}{dt}\|\rho_1 - \rho_2\|_{L^2}^2 &= -\intt (\rho_1 - \rho_2) \nabla\cdot [(\rho_1 - \rho_2)u_1 + \rho_2 (u_1-u_2)]\,dx\\
&= -\frac12\intt (\nabla\cdot u_1) (\rho_1 -\rho_2)^2\,dx  -\intt (\rho_1 - \rho_2) \nabla\rho_2 \cdot (u_1 - u_2)\,dx\\
&\quad - \intt (\rho_1 -\rho_2) \rho_2 \nabla\cdot (u_1 -u_2)\,dx\\
&\le C\|\nabla u_1\|_{L^\infty}\|\rho_1 -\rho_2\|_{L^2}^2 + \|\nabla\rho_2 \|_{L^3}\|\rho_1 - \rho_2\|_{L^2}\|u_1 - u_2\|_{L^6}\\
&\quad + C\|\rho_2\|_{L^\infty}\|\rho_1 -\rho_2\|_{L^2}\|\nabla (u_1 - u_2)\|_{L^2}\\
&\le C\|\rho_1 - \rho_2\|_{L^2}^2 + C\|u_1 - u_2\|_{H^1}^2,
\end{aligned}
\eq

Hence, we gather \eqref{eq:udiff_est_comb}, \eqref{eq:vdiff} and \eqref{eq:rhodiff} to have

\[\begin{aligned}
\begin{aligned}
&\frac{d}{dt}\lt( \|\rho_1 - \rho_2\|_{L^2}^2 + \|u_1 - u_2\|_{H^1}^2 + \frac{C_1}{\delta}\|v_1 - v_2\|_{L^2}^2\rt)\\
&\quad\le C\lt( \|\rho_1 - \rho_2\|_{L^2}^2 + \|u_1 - u_2\|_{H^1}^2 + \|v_1 - v_2\|_{L^2}^2\rt).
\end{aligned}
\end{aligned}\]
Thus, we use Gr\"onwall's lemma to get
\[
\begin{aligned}
&\lt( \|\rho_1 - \rho_2\|_{L^2}^2 + \|u_1 - u_2\|_{H^1}^2 + \|v_1 - v_2\|_{L^2}^2\rt)\\
&\quad\le C\lt( \|\rho_{1,0} - \rho_{2,0}\|_{L^2}^2 + \|u_{1,0} - u_{2,0}\|_{H^1}^2 + \|v_{1,0} - v_{2,0}\|_{L^2}^2\rt)
\end{aligned}
\]
and this implies the uniqueness.

\subsection{Local-in-time existence theory}
For the local-in-time existence, we decouple, regularize the system and apply Schauder's fixed point argument as follows:

We consider the Banach space $X_T = \calC([0,T]; H^2(\mathbb{T}^3)) \times \calC([0,T]; H^2(\mathbb{T}^3))$, where $T$ will be determined later. We set $R := 4(\|u_0\|_{H^3} + \|v_0\|_{H^3})$ and
\[
\mathcal{B}_R := \{(\bar u, \bar v) \in X_T \ : \ \|\bar u\|_{\calC([0,T];H^2(\T^3))} + \|\bar v\|_{\calC([0,T];H^2(\T^3))} \le R \}.
\]
For $(\bar u,\bar v)\in\mathcal{B}_R$, we define $\mathcal{T}(\bar u,\bar v):=(u,v)$ through the following decoupled, regularized system:

\begin{equation}\label{eq:decp_sys}
\begin{aligned}
&\pa_t \rho + \nabla\cdot (\rho (\bar u)_\e) = 0,\\
&\pa_t u + ((\bar u)_\e\cdot \nabla)u = -(u-(\bar v)_\e),\\
&\pa_t v + (v \cdot \nabla) v  + \nabla p -\nabla\cdot (G[|Dv|^2]Dv) = \rho(u-(\bar v)_\e),\\
&\nabla\cdot v= 0.
\end{aligned}
\end{equation}

subject to initial data:
\[
(\rho(0,x), u(0,x), v(0,x)) =: (\rho_0(x), u_0(x), v_0(x)), \quad x \in \T^3,
\]
where $(f)_\e := \theta_\e \star f$ denotes the mollification using the standard mollifier, i.e. $\theta_\e(x) := \sum_{k\in\mathbb{Z}^3}\frac{1}{\e^3}\tilde{\theta}\lt(\frac{x+k}{\e}\rt)$ where $\tilde{\theta}$ satisfies
\[
0 \leq \tilde{\theta} \in \mathcal{C}^\infty_0(\R^3), \quad \mbox{supp}_x\tilde{\theta} \subseteq B(0,1), \quad \int_{\R^3} \tilde{\theta}(x)\,dx = 1.
\]
 On the one hand,  the existence of solutions $(\rho, u)$ to \eqref{eq:decp_sys}$_1$-\eqref{eq:decp_sys}$_2$ follows from the classical theory of transport equations (see \cite{DL89} for example). On the other hand, adapting the local existence argument of \cite{KKK19} to the periodic domain with forcing $\rho(u-(\bar v)_\e)$ gives a solution $v$ to \eqref{eq:decp_sys}$_3$-\eqref{eq:decp_sys}$_4$.

First, we show that the operator $\mathcal{T}$ is well-defined, continuous and compact to use Schauder's fixed point theorem.
Similar estimates to Lemmas \ref{lem:l2est_v}-\ref{lem:l2est_h} give

\[
\frac{d}{dt}\|\rho\|_{H^2} \le C\|(\bar u)_\e \|_{H^3} \|\rho \|_{H^2},
\]
\[
\frac{d}{dt}\|u\|_{H^3}^2 + \|u\|_{H^3}^2 \le C\|(\bar u)_\e\|_{H^3} \|u\|_{H^3}^2 + \|(\bar v)_\e\|_{H^3}^2,
\]
and
\[\begin{aligned}
\frac{d}{dt}\|v\|_{H^3}^2 + 2c_0 \|\nabla v\|_{H^3}^2 &\le C\|v\|_{H^3}^3\\
&\quad + C\left(\max_{0\le s\le C\|v\|_{H^3}^2}G[s]\right)(1+ \|v\|_{H^3})\|v\|_{H^3}^{3/2}\|\nabla v\|_{H^3}^{3/2}\\
&\quad + C\|\rho\|_{H^2}(\|u\|_{H^2} + \|(\bar v)_\e\|_{H^2})(\|v\|_{L^2}+\|\nabla v\|_{H^3}).
\end{aligned}\]
Then we choose $T=T_\e$ sufficiently small, depending only on $\|\rho_0\|_{H^2}$, $\|u_0\|_{H^3}$, $\|v_0\|_{H^3}$, $R$ and $\e$, to have
\bq\label{eq:bd_rhou}
\begin{aligned}
&\sup_{0 \le t \le T_\e} \|\rho(t)\|_{H^2} \le 2\|\rho_0\|_{H^2},
\\
&\sup_{0 \le t \le T_\e} \|u(t)\|_{H^3} \le 2(\|u_0\|_{H^3}+\|v_0\|_{H^3}),
\end{aligned}
\eq
and
\bq\label{eq:bd_v}
\sup_{0 \le t \le T_\e} \|v(t)\|_{H^3}^2 + c_0\int_0^{T_\e} \|\nabla v (s)\|_{H^3}^2\,ds \le 4(\|u_0\|_{H^3}+\|v_0\|_{H^3})^2,
\eq

which implies the operator $\mathcal{T} : \mathcal{B}_R \to X_{T_\e}$ is well-defined. For the continuity, it follows from the uniform bounds \eqref{eq:bd_rhou} and \eqref{eq:bd_v} combined with $L^2$ stability estimates above. For the compactness, note that the uniform bounds from \eqref{eq:bd_rhou} and \eqref{eq:bd_v} also imply
\[
\|\pa_t u\|_{L^\infty([0,T_\e];H^2(\T^3))}, \quad \|\pa_t v\|_{L^2([0,T_\e]; H^2(\T^3))} <\infty.
\]
Hence, if we have a sequence $(\bar u_n, \bar v_n)$ in $\mathcal{B}_R$ and write the corresponding solution to \eqref{eq:decp_sys} as $(\rho_n, u_n, v_n)$, we can use the Aubin-Lions lemma to find that $(u_n, v_n)$ is relatively compact in $X_{T_\e}$, which shows the compactness of the operator.
Then we choose $T_\e$ smaller if necessary, to show that $\mathcal{T}$ maps $\mathcal{B}_R$ into $\mathcal{B}_R$. Since $\mathcal{T}$ is continuous and compact, it has a fixed point in $\mathcal{B}_R$ by Schauder's fixed point theorem, which yields a solution to the following regularized system on the time interval $[0,T_\e]$:

\begin{equation}\label{eq:reg_sys}
\begin{aligned}
&\pa_t \rho + \nabla\cdot (\rho ( u)_\e) = 0,\\
&\pa_t u + ((u)_\e\cdot \nabla)u = -(u-(v)_\e),\\
&\pa_t v + (v \cdot \nabla) v  + \nabla p -\nabla\cdot (G[|Dv|^2]Dv) = \rho(u-(v)_\e),\\
&\nabla\cdot v= 0.
\end{aligned}
\end{equation}

For \eqref{eq:reg_sys},  we again follow the estimates in Lemmas \ref{lem:l2est_v}-\ref{lem:l2est_h} to have

\[
\frac{d}{dt}\|\rho\|_{H^2} \le C\|(u)_\e \|_{H^3} \|\rho \|_{H^2} \le C\|u\|_{H^3} \|\rho\|_{H^2},
\]
\[
\frac{d}{dt}\|u\|_{H^3}^2 + \|u\|_{H^3}^2 \le C\|( u)_\e\|_{H^3} \|u\|_{H^3}^2 + \|(v)_\e\|_{H^3}^2 \le C\|u\|_{H^3}^3 + \|v\|_{H^3}^2,
\]
and
\[\begin{aligned}
\frac{d}{dt}\|v\|_{H^3}^2 + 2c_0 \|\nabla v\|_{H^3}^2 &\le C\|v\|_{H^3}^3\\
&\quad + C\left(\max_{0\le s\le C\|v\|_{H^3}^2}G[s]\right)(1+ \|v\|_{H^3})\|v\|_{H^3}^{3/2}\|\nabla v\|_{H^3}^{3/2}\\
&\quad + C\|\rho\|_{H^2}(\|u\|_{H^2} + \|( v)_\e\|_{H^2})(\|v\|_{L^2}+\|\nabla v\|_{H^3})\\
&\le C\|v\|_{H^3}^3 + C\|\rho\|_{H^2}(1+\|\rho\|_{H^2})(\|u\|_{H^2}^2 + \|v\|_{H^2}^2)+ \frac{c_0}{2}\|\nabla v\|_{H^3}^2\\
&\quad + C\left(\max_{0\le s\le C\|v\|_{H^3}^2}G[s]\right)(1+ \|v\|_{H^3})\|v\|_{H^3}^{3/2}\|\nabla v\|_{H^3}^{3/2},
\end{aligned}\]
and the right-hand sides are now independent of $\e$.  Hence, by using the standard continuation argument, we can construct a local-in-time solution $(\rho, u,v)$ to \eqref{eq:reg_sys} on the time interval $[0,T^*]$ satisfying
\[
\begin{aligned}
&\sup_{0 \le t \le T^*} \|\rho(t)\|_{H^2} \le 3\|\rho_0\|_{H^2},
\\
&\sup_{0 \le t \le T^*} \|u(t)\|_{H^3} \le 3(\|u_0\|_{H^3}+\|v_0\|_{H^3}),
\end{aligned}
\]
and
\[
\sup_{0 \le t \le T^*} \|v(t)\|_{H^3}^2 + c_0\int_0^{T^*} \|\nabla v (s)\|_{H^3}^2\,ds \le 9(\|u_0\|_{H^3}+\|v_0\|_{H^3})^2,
\]
where the choice of $T^*$ depends on $\|\rho_0\|_{H^2}$, $\|u_0\|_{H^3}$, and $\|v_0\|_{H^3}$, but independent of $\e$.  Then one uses these uniform bounds and Aubin-Lions lemma to establish a solution $(\rho, u, v)$ to \eqref{main_sys} with the desired regularity. The time continuity of $(\rho,u)$ in $H^2\times H^3$ follows from the transport energy estimates, and that of $v$ in $H^3$ follows from $v\in L^2(0,T^*;H^4)$ and $\pa_t v\in L^2(0,T^*;H^2)$. The positivity of $\rho$ follows from the characteristic formula.

Now, let $0<\e_0 <\delta_0$ be given and assume that
\[
\|\rho_0 -\rho_c(0)\|_{H^2} + \|u_0 - u_c(0)\|_{H^3} + \|v_0 - v_c(0)\|_{H^3} + |u_c(0)-v_c(0)| < \e_0.
\]
Here, we also let $\rho_c(0) =1$ and hence, $\rho_c(t) \equiv 1$ for any $t \in [0,T^*]$. Rather than using the energy inequality, we directly estimate $u-u_c$ and $v-v_c$ as follows:
\[\begin{aligned}
\frac12\frac{d}{dt}\|u-u_c\|_{L^2}^2 &= \frac12\intt (\nabla\cdot u) |u-u_c|^2\,dx - \intt (u-u_c)\cdot (u-v)\,dx,
\end{aligned}\]
and
\[\begin{aligned}
\frac12\frac{d}{dt}\|v-v_c\|_{L^2}^2  + c_0\|\nabla v\|_{L^2}^2 \le -\intt (v-v_c)\cdot (v-u)\,dx - \intt (\rho-1)(v-v_c)\cdot (v-u)\,dx.
\end{aligned}\]
These give
\[\begin{aligned}
\frac12&\frac{d}{dt}\lt(\|u-u_c\|_{L^2}^2 + \|v-v_c\|_{L^2}^2 \rt) + c_0\|\nabla v\|_{L^2}^2 + \|u-v\|_{L^2}^2\\
&\le \|\nabla u\|_{L^\infty}\|u-u_c\|_{L^2}^2 + |u_c - v_c|^2 + \|\rho-1\|_{L^\infty}\|v-v_c\|_{L^2}\|u-v\|_{L^2}\\
&\le C\|u-u_c\|_{H^3}^3 + |u_c - v_c|^2 + C\|\rho-1\|_{H^2}^2\|v-v_c\|_{L^2}^2 + \frac12\|u-v\|_{L^2}^2
\end{aligned}\]
We use $\intt (\rho-1)\,dx = 0$ to have
\[
\begin{aligned}
\frac12\frac{d}{dt}|u_c - v_c|^2 &= (u_c - v_c)\cdot \frac{d}{dt}(u_c - v_c)\\
&= -\intt (u_c -v_c) \cdot [(u\cdot \nabla) u]\,dx - 2|u_c - v_c|^2  -\intt (\rho-1)(u_c- v_c)(v-u)\,dx\\
&\le |u_c - v_c| \|u-u_c\|_{L^2}\|\nabla u\|_{L^2} - 2|u_c - v_c|^2  + \|\rho-1\|_{L^2} |u_c - v_c| \|v-u\|_{L^2}\\
&\le C|u_c - v_c| \|u-u_c\|_{H^3}^2 - 2|u_c - v_c|^2 + \|\rho-1\|_{H^2}^2 |u_c -v_c|^2 + \frac12\|v-u\|_{L^2}^2.
\end{aligned}
\]
Collecting the above estimates, we attain

\[
\begin{aligned}
\frac{d}{dt}&\lt( \|u-u_c\|_{L^2}^2 + \|v-v_c\|_{L^2}^2 + |u_c - v_c|^2\rt)  + 2c_0 \|\nabla v\|_{L^2}^2\\
&\le C\|u-u_c\|_{H^3}^3 +  C\|\rho-1\|_{H^2}^2\lt( \|v-v_c\|_{L^2}^2 + |u_c - v_c|^2\rt) + C|u_c - v_c| \|u-u_c\|_{H^3}^2.
\end{aligned}
\]
Then we combine the above inequality with Lemmas \ref{lem:l2est_v}, \ref{lem:l2est_u}, and \ref{lem:l2est_h} to obtain, 
\[\begin{aligned}
\frac{d}{dt}&\lt(\|\rho-\rho_c\|_{H^2}^2 + \|u-u_c\|_{H^3}^2 + \|v-v_c\|_{H^3}^2 + |u_c - v_c|^2 \rt) + c_0 \|\nabla v\|_{H^3}^2\\
&\le F\lt( \|\rho-\rho_c\|_{H^2}^2 + \|u-u_c\|_{H^3}^2 + \|v-v_c\|_{H^3}^2 + |u_c - v_c|^2\rt),
\end{aligned}\]
for some nonnegative continuous function $F : \R \to [0,\infty) $ satisfying $F(0) = 0$. Then reducing $T^*$ if necessary, we can get 
\[
\sup_{t \in [0,T^*]} \lt( \|\rho(t, \cdot) -\rho_c(t)\|_{H^2} + \|u(t,\cdot) - u_c(t)\|_{H^3} + \|v(t,\cdot) - v_c(t)\|_{H^3} + |u_c(t) - v_c(t)|\rt)<\delta_0.
\]
This completes the proof.

\end{document}